\documentclass[11pt,a4paper]{amsart}
\usepackage[T1]{fontenc}
\usepackage{lmodern}
\usepackage{amsmath,amssymb,amsthm,mathtools}
\usepackage{hyperref}
\hypersetup{colorlinks=true,linkcolor=blue,citecolor=blue,urlcolor=blue}
\usepackage{microtype}

\newtheorem{theorem}{Theorem}[section]
\newtheorem{proposition}[theorem]{Proposition}
\newtheorem{corollary}[theorem]{Corollary}

\theoremstyle{remark}
\newtheorem{remark}[theorem]{Remark}
\newtheorem{example}[theorem]{Example}
\newcommand{\R}{\mathbb{R}} 
\newcommand{\Oc}{\mathcal O_{c}}
\newcommand{\Oo}{\mathcal O_{o}}
\newcommand{\Bc}{\mathcal B_{c}} 
\newcommand{\Bo}{\mathcal B_{o}}
\newcommand{\rank}{\operatorname{rank}}
\newcommand{\col}{\operatorname{col}}
\newcommand{\epsb}{\bar\varepsilon}
\newcommand{\zm}{\underline z^{m}}
\newcommand{\zn}{\underline z^{n}}
\newcommand{\zp}{\underline{\bar z}^{*}}

\title[Multirate predictors: rank and conditioning]
{Multirate Data-Driven Predictive Control: Rank Obstructions, Identifiability, and Noise Conditioning}
\author{Robert Vrabel}
\address{Institute of Applied Informatics, Automation and Mechatronics, Slovak University of Technology in Bratislava, Bottova 25, 917 24 Trnava, Slovakia}
\email{robert.vrabel@stuba.sk}

\begin{document} 
\begin{abstract}
We develop rank, identifiability, and noise-conditioning results for compact input--output predictors in multirate data-driven predictive control. The inverse of a square zero-input response matrix used in an existing compact-lifting construction is shown to be undefined under its observability and minimal-polynomial assumptions. We derive exact replacements using rank-complete data, characterize all predictor extensions, and establish the minimum-norm and minimax Euclidean noise properties of the canonical Moore--Penrose extension. A joint predictor is identified directly from mixed experiments with arbitrary initial states. Deterministic perturbation bounds quantify reconstruction errors through informative singular-value margins, while a positive-definite reduced weighting preserves the nominal quadratic cost. Sampling-rate observability and event-triggering analyses identify additional conditions needed for reconstruction and stability guarantees.
\end{abstract}
\maketitle

\section{Setting and purpose}
Consider the continuous-time linear system
\begin{equation}\label{eq:ct}
  \dot x=Ax+Bu,\qquad y=Cx,
\end{equation}
with state dimension $n$, input dimension $m$, and output dimension $p\ge 1$.
The input is held over intervals of length $T$, whereas the output is sampled $h$ times per input interval. Write 
\begin{align*}
 A_d&=e^{AT}, & B_d&=\int_0^T e^{A\tau}B\,d\tau,\\
 A_f&=e^{AT/h}, & B_f&=\int_0^{T/h}e^{A\tau}B\,d\tau.
\end{align*} 
Behavioral input--output representations based on persistently exciting data underlie Data-Enabled Predictive Control (DeePC)~\cite{Willems2005,Coulson2019}.
Subspace predictive control, instrumental-variable data representations, and connections between direct and indirect data-driven methods provide related predictor-based formulations~\cite{Favoreel1999,vanWingerden2022,Dorfler2023,Dinkla2026}.
In a multirate setting, an input may remain constant over a slow interval while output measurements within that interval permit a compact representation of the remaining samples. A compact multirate output-lifting construction in Yang et al.~\cite{Yang2025} provides a concrete setting for examining the role of data rank and observability. Its standing assumptions include controllability of $(A_d,B_d)$ and observability of $(A_f,C)$. Let $r$ be the degree of the minimal polynomial of $A_f$. We consider nontrivial compact lifting with $h>r+1$. Put $s=(r+1)p$ and $q=(h-r-1)p>0$; the first $r+1$ fast output samples form a retained block in $\R^s$, and the remaining samples form an omitted block in $\R^q$. The resulting predictor is an exact map on a low-dimensional realizable subspace, whose ambient-space extensions have different noise and conditioning properties.

The main contributions are as follows.
\begin{enumerate}
\item We establish a universal rank obstruction to an inverse-based zero-input compact-lifting construction and a rank-correct exact replacement.
\item We characterize all linear extensions of the realizable output map and prove the minimum-norm and minimax Euclidean noise properties of its canonical pseudoinverse extension.
\item We derive a direct predictor from mixed input--output experiments with arbitrary initial states and a measurable joint-data rank condition.
\item We prove deterministic perturbation estimates for identified lifting and input-response matrices, exposing the singular-value margins controlling sensitivity.
\item We distinguish slow- and fast-rate observability and establish precise limitations of a finite-horizon event trigger and a subsequent stability-proof estimate.
\end{enumerate}

The results distinguish three levels of analysis: exact prediction on the realizable subspace, extension and identification from finite experimental data, and stability of the resulting closed-loop implementation. Equation numbers attributed to~\cite{Yang2025} refer to that publication.
 
\section{An unavoidable rank obstruction in inverse-based lifting}
For a zero-input fast-sampled trajectory starting from $x$, define its retained and omitted output blocks as
\begin{align}
 z_c&=\Oc x, & z_o&=\Oo x,\label{eq:obs}\\
 \Oc&=\col(C,CA_f,\ldots,CA_f^r),\notag\\
 \Oo&=\col(CA_f^{r+1},\ldots,CA_f^{h-1}).\notag
\end{align}
The offline procedure in the proof of \cite[Theorem 1]{Yang2025} selects $s$ such responses, puts them into the square matrix
\[
 \Phi_x=[z_{c,1},\ldots,z_{c,s}]\in\R^{s\times s},
\]
and declares $\Phi_x$ nonsingular in order to compute $G_x=\Phi_y\Phi_x^{-1}$, where $\Phi_y=[z_{o,1},\ldots,z_{o,s}]$.

\begin{theorem}\label{thm:singular}
Suppose $(A_f,C)$ is observable and the minimal polynomial of $A_f$ has degree $r$. For every choice of initial states $x_1,\ldots,x_s$, where $s=(r+1)p$, one has
\begin{equation}\label{eq:strongrank}
 \rank\Phi_x\le n\le rp<(r+1)p=s.
\end{equation}
Consequently, the square zero-input response matrix $\Phi_x$ is singular for every system satisfying these assumptions, irrespective of the offline excitation.
\end{theorem}
\begin{proof}
Since the minimal polynomial has degree $r$, every power $A_f^k$ with $k\ge r$ is a linear combination of $I,A_f,\ldots,A_f^{r-1}$. Hence all later block rows of the observability matrix are linear combinations of the first $r$ block rows. Observability therefore implies
\[
 \rank\col(C,CA_f,\ldots,CA_f^{r-1})=n,
\]
so $n\le rp$. With $X=[x_1,\ldots,x_s]\in\R^{n\times s}$, one has $\Phi_x=\Oc X$ by \eqref{eq:obs}; thus $\rank\Phi_x\le n$. Finally $rp<(r+1)p$ because $p\ge1$.
\end{proof}

The nominal compact-lifting relation exists on the realizable retained-output subspace. The inverse-based identification $G_x=\Phi_y\Phi_x^{-1}$ in \cite[Theorem~1]{Yang2025} is undefined for every admissible observable realization. A constructive extension is given next.

\begin{example}\label{ex:two}
Take
\[
 A=\begin{bmatrix}-1&0\\0&-2\end{bmatrix},\qquad
 B=\begin{bmatrix}1\\1\end{bmatrix},\qquad
 C=\begin{bmatrix}1&1\end{bmatrix},\qquad T=1,\quad h=4.
\]
Here $A_d=\operatorname{diag}(e^{-1},e^{-2})$ and
$B_d=\col(1-e^{-1},(1-e^{-2})/2)$, so
\[
 \det[B_d,A_dB_d]
 =\frac{(1-e^{-1})(1-e^{-2})}{2}(e^{-2}-e^{-1})\ne0.
\]
Also $A_f=\operatorname{diag}(e^{-1/4},e^{-1/2})$ and
\[
 \det\col(C,CA_f)=e^{-1/2}-e^{-1/4}\ne0.
\]
Thus the controllability and observability conditions hold, $r=2$, and $r+1=3<h$. A scalar offline input persistently exciting of the required order $L+n+1$ can be chosen for sufficiently large data length $N$; this does not change the bound $\rank\Phi_x\le2<3$.
\end{example}

\section{A rank-correct compact lifting}
The same observability argument that establishes Theorem~\ref{thm:singular} provides a constructive alternative. Define the zero-order-hold output-response blocks
\[
 b_0=0_{p\times m},\qquad
 b_j=C\sum_{\nu=0}^{j-1}A_f^{\nu}B_f\in\R^{p\times m},
 \quad j=1,\ldots,h-1,
\]
and set
\[
 \Bc=\col(b_0,\ldots,b_r)\in\R^{s\times m},\qquad
 \Bo=\col(b_{r+1},\ldots,b_{h-1})\in\R^{q\times m}.
\]
For arbitrary $x$ and constant input $u$ over a slow interval,
\begin{equation}\label{eq:inputoutput}
 z_c=\Oc x+\Bc u,\qquad z_o=\Oo x+\Bo u.
\end{equation}

\begin{proposition}\label{prop:correct}
Assume the hypotheses of Theorem~\ref{thm:singular}. Let $H$ be a real matrix of size $(h-r-1)p\times(r+1)p$ satisfying
\begin{equation}\label{eq:intertwine}
 H\Oc=\Oo.
\end{equation}
Then every nominal trajectory satisfies
\begin{equation}\label{eq:extension}
 z_o=Hz_c+Ku,\qquad K:=\Bo-H\Bc.
\end{equation}
One such choice is $H=\Oo\Oc^{\dagger}$. Alternatively, form zero-input data matrices $\Phi_x=\Oc X$ and $\Phi_y=\Oo X$, with the same initial-state matrix $X$. If
\begin{equation}\label{eq:rankdata}
 \rank\Phi_x=n,
\end{equation}
then the entirely data-defined choice $H=\Phi_y\Phi_x^{\dagger}$ satisfies \eqref{eq:intertwine} and hence \eqref{eq:extension}.
\end{proposition}
\begin{proof}
Observability and the minimal-polynomial relation give $\rank\Oc=n$, so $\Oc^{\dagger}\Oc=I_n$ and the first choice satisfies \eqref{eq:intertwine}. If \eqref{eq:rankdata} holds, then $\rank X=n$ because $\Oc$ has full column rank. Let $X^R$ be a right inverse of $X$. Since $\Phi_y=\Oo X$, the Moore--Penrose identity and $\operatorname{row}\Phi_y\subseteq\operatorname{row}\Phi_x$ yield
\[
 H\Phi_x=\Phi_y\Phi_x^{\dagger}\Phi_x=\Phi_y.
\]
Consequently $H\Oc=H\Oc XX^R=\Oo XX^R=\Oo$. Substituting $\Oc x=z_c-\Bc u$ into $z_o=\Oo x+\Bo u$ proves \eqref{eq:extension}.
\end{proof}

The matrix $K$ can also be identified from exact input--output responses with arbitrary initial states. For $N_u$ experiments, let $X_U=[x_1,\ldots,x_{N_u}]$ contain the unknown initial states, and let $Z_c$, $Z_o$, and $U=[u_1,\ldots,u_{N_u}]\in\R^{m\times N_u}$ collect the retained outputs, omitted outputs, and constant inputs. From \eqref{eq:inputoutput} and \eqref{eq:intertwine},
\begin{equation}\label{eq:Kdata}
 Z_o-HZ_c=(\Oo-H\Oc)X_U+(\Bo-H\Bc)U=KU.
\end{equation}
Consequently, if $\rank U=m$, then
\begin{equation}\label{eq:Kpseudoinverse}
 K=(Z_o-HZ_c)U^\dagger.
\end{equation}
In particular, the unit inputs $u_j=e_j$ recover the columns of $K$ individually. Both $H$ and $K$ therefore admit data-based identification with the stated rank conditions and exact nominal responses, without prescribed initial states for the input experiments.

Let $Q\in\R^{(s+q)\times(s+q)}$ and $R\in\R^{m\times m}$ be symmetric positive-definite full-lifted output and input weights. Define
\[
 M_H=\begin{bmatrix}I_s&0_{s\times m}\\ H&K\end{bmatrix}
 \in\R^{(s+q)\times(s+m)},
 \qquad
 G_c=M_H^\top Q M_H+\operatorname{diag}(0_{s\times s},R).
\]
This $G_c$ is positive definite: if the quadratic form vanishes, the $R$ term gives $u=0$, and the upper identity block of $M_H$ then gives $z_c=0$. Furthermore, for every nominal trajectory,
\begin{equation}\label{eq:cost}
 \begin{bmatrix}z_c\\u\end{bmatrix}^{\!\top}G_c
 \begin{bmatrix}z_c\\u\end{bmatrix}
 =\begin{bmatrix}z_c\\z_o\end{bmatrix}^{\!\top}Q
 \begin{bmatrix}z_c\\z_o\end{bmatrix}+u^\top Ru.
\end{equation}
The construction establishes an exact nominal compact-lifting identity using the rank-$n$ data and its Moore--Penrose inverse.

\begin{remark}
The equality \eqref{eq:cost} holds for nominal trajectories. For measured blocks $z_c+\varepsilon_c$ and $z_o+\varepsilon_o$, the extension residual is
\begin{equation}\label{eq:residual}
 (z_o+\varepsilon_o)-H(z_c+\varepsilon_c)-Ku
 =\varepsilon_o-H\varepsilon_c.
\end{equation}
Under the componentwise measurement-noise bound $\|\varepsilon\|_\infty\le\bar\varepsilon$ of \cite[Assumption~2]{Yang2025}, the Euclidean residual satisfies
\begin{equation}\label{eq:noisecalc}
 \|\varepsilon_o-H\varepsilon_c\|_2
 \le \bar\varepsilon\bigl(\sqrt{(h-r-1)p}
  +\|H\|_2\sqrt{(r+1)p}\bigr).
\end{equation}
This estimate describes measurement-noise propagation for matrices $H$ and $K$ identified from exact nominal responses. Identification errors in these matrices introduce additional perturbation terms, governed by their own matrix-estimation bounds.
\end{remark}

\section{Canonical extensions and optimal Euclidean noise amplification}
An exact lifting map is prescribed only on the $n$-dimensional realizable output subspace. Its action on the orthogonal complement is not determined by the system's nominal trajectories. The following characterization makes this freedom explicit.

\begin{theorem}\label{thm:canonical}
Let $S\in\R^{s\times n}$ have full column rank and let $T_o\in\R^{q\times n}$, with $q\ge1$. Put
\[
 P=SS^\dagger,\qquad H_0=T_oS^\dagger.
\]
Every matrix $H\in\R^{q\times s}$ satisfying $HS=T_o$ has the representation
\begin{equation}\label{eq:allH}
 H=H_0+F(I_s-P),\qquad F\in\R^{q\times s}.
\end{equation}
Moreover,
\begin{align}
 HH^\top&=H_0H_0^\top+F(I_s-P)F^\top,\label{eq:HH}\\
 \|H\|_F^2&=\|H_0\|_F^2+\|F(I_s-P)\|_F^2.\label{eq:HF}
\end{align}
Consequently, $H_0$ is the unique exact extension of smallest Frobenius norm and is an exact extension of smallest spectral norm:
\begin{equation}\label{eq:specoptimal}
 \|H\|_2\ge \|H_0\|_2.
\end{equation}
For additive retained- and omitted-output errors satisfying
$\|[\varepsilon_o^\top,\varepsilon_c^\top]^\top\|_2\le\eta$,
\begin{equation}\label{eq:minimax}
 \sup_{\left\|\left[\begin{smallmatrix}\varepsilon_o\\\varepsilon_c\end{smallmatrix}\right]\right\|_2\le\eta}
 \|\varepsilon_o-H\varepsilon_c\|_2
 =\eta\sqrt{1+\|H\|_2^2}.
\end{equation}
Thus $H_0$ also minimizes this worst-case Euclidean amplification among all exact extensions.
\end{theorem}
\begin{proof}
Full column rank gives $S^\dagger S=I_n$, so $H_0S=T_o$. If $HS=T_o$, then $(H-H_0)S=0$, hence $(H-H_0)P=0$ and $H-H_0=(H-H_0)(I_s-P)$. This gives \eqref{eq:allH}. The identities $H_0(I_s-P)=0$ and $(I_s-P)^2=I_s-P$ eliminate cross terms in $HH^\top$ and in the Frobenius inner product, giving \eqref{eq:HH}--\eqref{eq:HF}. Positive semidefiniteness of the second term in \eqref{eq:HH} yields \eqref{eq:specoptimal}. Finally, the residual operator on the joint noise vector is $[I_q,-H]$, with
\[
 [I_q,-H][I_q,-H]^\top=I_q+HH^\top.
\]
Its spectral norm equals $\sqrt{1+\|H\|_2^2}$, proving \eqref{eq:minimax}.
\end{proof}

\begin{remark}
The projector $P$ separates a retained-output vector $v$ into its physically realizable component $Pv$ and an orthogonal component $(I_s-P)v$. Every exact extension has the same action on $Pv$. Its freedom is confined to $(I_s-P)v$ through the term $F(I_s-P)$ in \eqref{eq:allH}. The canonical choice $H_0$ satisfies $H_0v=H_0Pv$, so an off-subspace perturbation of the retained outputs has no direct effect on the reconstructed omitted outputs. This geometric property underlies its minimum-norm and Euclidean minimax characterizations.
\end{remark}

For the compact multirate predictor, set $S=\Oc$ and $T_o=\Oo$. In particular,
\begin{equation}\label{eq:canonH}
 H_0=\Oo\Oc^\dagger.
\end{equation}
The data-based construction of Proposition~\ref{prop:correct} recovers this \emph{same canonical matrix} whenever $\rank\Phi_x=n$.

\begin{corollary}\label{cor:dataind}
If $\Phi_x=\Oc X$ and $\Phi_y=\Oo X$, with $\rank\Phi_x=n$, then
\begin{equation}\label{eq:datainvariance}
 \Phi_y\Phi_x^\dagger=\Oo\Oc^\dagger.
\end{equation}
Hence every sufficiently rank-rich exact zero-input data set yields the same minimum-Frobenius-norm extension, independently of the number and selection of its response columns.
\end{corollary}
\begin{proof}
Proposition~\ref{prop:correct} implies that $H_d=\Phi_y\Phi_x^\dagger$ satisfies $H_d\Oc=\Oo$. Since $\operatorname{im}\Phi_x=\operatorname{im}\Oc$, the pseudoinverse property gives $H_d=H_dP$, with $P=\Oc\Oc^\dagger$. Representation \eqref{eq:allH} then forces $H_d=\Oo\Oc^\dagger$.
\end{proof}

\begin{remark}
The spectral-norm minimizer need not be unique; uniqueness in Theorem~\ref{thm:canonical} concerns the Frobenius norm. The minimax statement assumes a joint Euclidean noise ball. A componentwise noise set can induce a different optimal extension. The canonical extension minimizes the uncertainty gain of the \emph{linear reconstruction map}; stability and robustness of a receding-horizon closed loop also depend on its optimization constraints and dynamics.
\end{remark}

\section{Direct identification from mixed input--output experiments}
Proposition~\ref{prop:correct} identifies $H$ from zero-input responses and $K$ from general input experiments. A single collection of mixed experiments can identify the entire predictor without recovering its components sequentially.

With the retained and omitted block dimensions $s$ and $q$ fixed above, define
\begin{equation}\label{eq:TF}
 \mathcal T=\begin{bmatrix}\Oc&\Bc\\0&I_m\end{bmatrix}
 \in\R^{(s+m)\times(n+m)},\qquad
 \mathcal W=\begin{bmatrix}\Oo&\Bo\end{bmatrix}\in\R^{q\times(n+m)}.
\end{equation}
For $N$ experiments with arbitrary initial states and constant input vectors, collect
\[
 Z_c=[z_{c,1},\ldots,z_{c,N}],\quad
 Z_o=[z_{o,1},\ldots,z_{o,N}],\quad
 U=[u_1,\ldots,u_N].
\]
The measured joint data matrix is
\[
 D=\begin{bmatrix}Z_c\\ U\end{bmatrix}\in\R^{(s+m)\times N}.
\]

\begin{theorem}\label{thm:joint}
Suppose $\rank\Oc=n$, the experiments are exact, and
\begin{equation}\label{eq:jointrank}
 \rank D=n+m.
\end{equation}
Then
\begin{equation}\label{eq:Pi}
 \Pi_0=Z_oD^\dagger=\mathcal W\mathcal T^\dagger
\end{equation}
is the unique minimum-Frobenius-norm extension of the exact predictor on $\operatorname{im}\mathcal T$ to a linear map from $\R^{s+m}$ into $\R^q$; it also minimizes the spectral norm among exact extensions. For every nominal state $x$ and input $u$,
\begin{equation}\label{eq:jointpredict}
 z_o=\Pi_0\begin{bmatrix}z_c\\u\end{bmatrix}.
\end{equation}
Writing $\Pi_0=[H_J,K_J]$, one obtains
\begin{equation}\label{eq:jointconstraints}
 H_J\Oc=\Oo,\qquad H_J\Bc+K_J=\Bo.
\end{equation}
No prescribed initialization of the experimental states is needed.
\end{theorem}
\begin{proof}
Let $X=[x_1,\ldots,x_N]$ and $V=\operatorname{col}(X,U)\in\R^{(n+m)\times N}$. By \eqref{eq:inputoutput},
\[
 D=\mathcal T V,\qquad Z_o=\mathcal W V.
\]
Because $\Oc$ has full column rank, $\mathcal T$ has full column rank $n+m$: its lower block first forces the input coordinate to zero, and its upper block then forces the state coordinate to zero. Thus \eqref{eq:jointrank} implies $\rank V=n+m$ and $\operatorname{im}D=\operatorname{im}\mathcal T$. With a right inverse $V^R$ of $V$,
\[
 \Pi_0 D=Z_oD^\dagger D=Z_o,
 \qquad \Pi_0\mathcal T=\Pi_0DV^R=Z_oV^R=\mathcal W,
\]
where $\operatorname{row}Z_o\subseteq\operatorname{row}D$ establishes the first identity. Consequently \eqref{eq:jointpredict} holds and block multiplication proves \eqref{eq:jointconstraints}. Since $D^\dagger$ annihilates $(\operatorname{im}D)^\perp$, the extension $\Pi_0$ is zero on $(\operatorname{im}\mathcal T)^\perp$. Therefore $\Pi_0=\Pi_0\mathcal T\mathcal T^\dagger=\mathcal W\mathcal T^\dagger$, proving the second identity in \eqref{eq:Pi}. The argument of Theorem~\ref{thm:canonical}, applied to $(S,T_o)=(\mathcal T,\mathcal W)$, gives the claimed unique Frobenius minimizer.
\end{proof}

The measurable rank condition \eqref{eq:jointrank} involves both state-response information and independent input directions. It is stronger than $\rank U=m$ alone, since the latter does not ensure that the experiments resolve every initial-state direction. It requires at least $n+m$ informative columns, without requiring the experimental initial states to be measured.

\begin{remark}
The data columns may be collected from consecutive slow input intervals in a single open-loop or closed-loop experiment. Each interval contributes its retained and omitted fast-output blocks and the corresponding held input. The construction requires the same system matrices across these intervals and the verifiable joint-rank condition \eqref{eq:jointrank}; in particular, feedback-generated inputs are admissible when the collected data attain that rank. No reset to a prescribed initial state is needed. The equality in \eqref{eq:Pi} also shows that any rank-complete exact data set identifies the same canonical joint predictor.
\end{remark}

\begin{remark}
The two-stage pair $[H_0,K]$ and the joint minimum-norm extension $\Pi_0=[H_J,K_J]$ agree on $\operatorname{im}\mathcal T$ but may differ elsewhere. Their canonical normalizations refer to different ambient spaces: $H_0$ minimizes the norm of the retained-output extension, while $\Pi_0$ minimizes the norm of the joint $[z_c^\top,u^\top]^\top$ extension. Both preserve the nominal relation and, with the corresponding $M_H$ or $M_{H_J}$, the quadratic identity \eqref{eq:cost}.
\end{remark}

\section{Noise conditioning of data-identified predictors}
Exact data rank and stable identification are distinct requirements. The following bounds build on standard singular-value perturbation principles~\cite{Wedin1972} and apply to compact output lifting and the joint predictor.

\begin{theorem}\label{thm:perturb}
Let $A_0\in\R^{d\times k}$ have full column rank, let $B_0\in\R^{q\times k}$, and define
\[
 H_0=B_0A_0^\dagger,\qquad P_0=A_0A_0^\dagger,\qquad
 \gamma=\sigma_{\min}(A_0)>0.
\]
Suppose the identified matrices are formed from perturbed experiments,
\[
 \widetilde A=A_0+E_c,\qquad
 \widetilde B=B_0+E_o,\qquad
 \widehat H=\widetilde B\widetilde A^\dagger,
\]
where $\rho_c=\|E_c\|_2<\gamma$ and $\rho_o=\|E_o\|_2$. Then $\widetilde A$ has full column rank and
\begin{equation}\label{eq:pertmain}
 \| (\widehat H-H_0)P_0\|_2
 \le \frac{1}{\gamma}\left[\rho_o+
 \frac{(\|B_0\|_2+\rho_o)\rho_c}{\gamma-\rho_c}\right].
\end{equation}
Moreover,
\begin{equation}\label{eq:pertoff}
 \|\widehat H(I_d-P_0)\|_2
 \le \frac{(\|B_0\|_2+\rho_o)\rho_c}{(\gamma-\rho_c)^2},
\end{equation}
so the sum of the right sides of \eqref{eq:pertmain} and \eqref{eq:pertoff} bounds $\|\widehat H-H_0\|_2$ on the entire ambient space.
\end{theorem}
\begin{proof}
Weyl's singular-value inequality yields
\[
 \sigma_{\min}(\widetilde A)\ge\gamma-\rho_c>0,
 \qquad
 \|\widehat H\|_2\le
 \frac{\|B_0\|_2+\rho_o}{\gamma-\rho_c}.
\]
Since $\widehat H\widetilde A=\widetilde B$ and $H_0A_0=B_0$,
\[
 (\widehat H-H_0)A_0=E_o-\widehat H E_c.
\]
Multiplication on the right by $A_0^\dagger$ gives
\[
 (\widehat H-H_0)P_0=(E_o-\widehat H E_c)A_0^\dagger,
\]
which proves \eqref{eq:pertmain}.
For the complementary component, let $\widetilde P=\widetilde A\widetilde A^\dagger$. We have
\[
 \|(I_d-P_0)\widetilde P\|_2
 =\|(I_d-P_0)E_c\widetilde A^\dagger\|_2
 \le\frac{\rho_c}{\gamma-\rho_c}.
\]
The orthogonal projectors $P_0$ and $\widetilde P$ have the same rank, so the principal-angle identity gives
$\|\widetilde P(I_d-P_0)\|_2=\|(I_d-P_0)\widetilde P\|_2$.
Hence
\[
 \|\widetilde A^\dagger(I_d-P_0)\|_2
 \le\frac{\rho_c}{(\gamma-\rho_c)^2}.
\]
As $H_0(I_d-P_0)=0$ and $\widehat H=\widetilde B\widetilde A^\dagger$, \eqref{eq:pertoff} follows. Finally split $\widehat H-H_0=(\widehat H-H_0)P_0+(\widehat H-H_0)(I_d-P_0)$.
\end{proof}

\begin{corollary}\label{cor:inputpert}
Suppose exact data satisfy $Z_o-HZ_c=KU$ with $\rank U=m$, while the output data used to identify $K$ are $\widetilde Z_c=Z_c+E_{c,K}$, $\widetilde Z_o=Z_o+E_{o,K}$, and $\widehat H=H+\Delta H$. If the input data $U$ are exact and
\[
 \widehat K=(\widetilde Z_o-\widehat H\widetilde Z_c)U^\dagger,
\]
then
\begin{equation}\label{eq:Kpert}
 \|\widehat K-K\|_2\le
 \frac{\|E_{o,K}\|_2+\|\widehat H\|_2\|E_{c,K}\|_2+
       \|\Delta H\|_2\|Z_c\|_2}{\sigma_{\min}(U)}.
\end{equation}
\end{corollary}
\begin{proof}
The defining identities and $UU^\dagger=I_m$ yield
\[
 \widehat K-K=
 (E_{o,K}-\widehat H E_{c,K}-\Delta H Z_c)U^\dagger.
\]
Taking spectral norms proves \eqref{eq:Kpert}.
\end{proof}

\begin{remark}
The singular-value margin $\gamma=\sigma_{\min}(A_0)$ measures the separation of the informative experimental directions from rank deficiency. Full column rank establishes exact identifiability, while the comparison $\rho_c<\gamma$ controls sensitivity to finite-data perturbations. As $\rho_c$ approaches $\gamma$, the factors $(\gamma-\rho_c)^{-1}$ and $(\gamma-\rho_c)^{-2}$ in \eqref{eq:pertmain}--\eqref{eq:pertoff} grow, signaling possible severe noise amplification. Theorem~\ref{thm:perturb} applies directly to $k$ independent columns, with $k=n$ for zero-input lifting and $k=n+m$ for the joint predictor. For larger noisy collections, the use of a rank-$k$ subspace model or truncated SVD calls for a corresponding perturbation analysis, since generic noise can increase the numerical rank of an untruncated data matrix.
\end{remark}

\begin{remark}
For a new measurement pair $(\widetilde z_c,\widetilde z_o)$ and input $u$, suppose the true predictor is $(H,K)$ and the learned matrices are $(\widehat H,\widehat K)=(H+\Delta H,K+\Delta K)$. With measurement errors $\varepsilon_c,\varepsilon_o$,
\begin{align*}
 &\widetilde z_o-\widehat H\widetilde z_c-\widehat Ku\\
 &\qquad=\varepsilon_o-\widehat H\varepsilon_c
  -\Delta H z_c-\Delta K u.
\end{align*}
Under a componentwise measurement-noise bound $\bar\varepsilon$, with $q_c=s$ and $q_o=q$, the residual satisfies
\begin{align}\label{eq:fullnoise}
 &\|\widetilde z_o-\widehat H\widetilde z_c-\widehat Ku\|_2
 \le \bar\varepsilon\bigl(\sqrt{q_o}+\|\widehat H\|_2\sqrt{q_c}\bigr)\notag\\
 &\hspace{12mm}+\|\Delta H\|_2\|z_c\|_2+\|\Delta K\|_2\|u\|_2.
\end{align}
Combining \eqref{eq:pertmain}--\eqref{eq:Kpert} with \eqref{eq:fullnoise} gives a deterministic error budget whenever the required experimental singular-value margins and signal bounds are available.
\end{remark}

\section{Slow-rate versus fast-rate observability}
In \cite[Definition 1]{Yang2025}, the observability index $\ell$ is defined using the slow-rate matrix $A_d$:
\begin{equation}\label{eq:slow}
 \ell_s:=\min\bigl\{j\ge1:\rank\col(C,CA_d,\ldots,CA_d^{j-1})=n\bigr\}.
\end{equation}
Yet \cite[Lemma 2]{Yang2025} reconstructs the initial state from fast samples and treats $\col(C,CA_f,\ldots,CA_f^{\ell_s-1})$ as full column rank. The assumed observability of $(A_f,C)$ does not ensure that the slow index exists.

\begin{example}\label{ex:rotation}
Set
\[
 A=\begin{bmatrix}0&-1\\1&0\end{bmatrix},\quad B=I_2,\quad
 C=\begin{bmatrix}1&0\end{bmatrix},\qquad T=\pi,\quad h=4.
\]
Then
\[
 A_d=-I_2,\qquad
 B_d=\int_0^\pi e^{A\tau}\,d\tau=
 \begin{bmatrix}0&-2\\2&0\end{bmatrix}.
\]
Thus $(A_d,B_d)$ is controllable, but $(A_d,C)$ is not observable. On the other hand, with the fast-sampling rotation angle $\varphi=\pi/4$,
\[
 A_f=\begin{bmatrix}\cos\varphi&-\sin\varphi\\\sin\varphi&\cos\varphi\end{bmatrix},\qquad
 \det\col(C,CA_f)=-\sin\varphi\ne0.
\]
Therefore $(A_f,C)$ is observable. Its minimal polynomial has degree $r=2$, so $r+1=3<h=4$ as required for nontrivial compact lifting. Nevertheless $\ell_s$ in \eqref{eq:slow} does not exist.
\end{example}

For state reconstruction from the available fast output samples, the relevant rank test must use the corresponding fast-sampled pair. The sampling-consistent reconstruction index is
\begin{equation}\label{eq:fast}
 \ell_f:=\min\bigl\{j\ge1:\rank\col(C,CA_f,\ldots,CA_f^{j-1})=n\bigr\}\le r.
\end{equation}
When a known input is held constant within each slow interval, $\ell_f$ successive fast outputs determine the initial state after subtracting the known input response. The number of slow input blocks covering these observations is $\lceil\ell_f/h\rceil$; under $r+1\le h$, this number equals one.

\section{Event well-posedness and the direction of the stability estimate}
Let $\epsb>0$ be a prescribed bound on measurement noise, even if the actual noise is zero. Use the notation of \cite{Yang2025} for the measured, nominal, and optimal predicted retained output sequences, $\zm$, $\zn$, and $\zp$, respectively. At an optimization instant, define the measured--predicted difference and the event threshold by
\begin{align}
 D_m(t)&=\|\zm_{[0,t-1]}\|^2-\|\zp_{[0,t-1]}\|^2,\notag\\
 \Theta(t)&=\delta_1\epsb+\delta_2\|\zp_{t-1}\|^2
 -2(r+1)t\epsb\|\zm_{[0,t-1]}\|,\label{eq:threshold}
\end{align}
where $\delta_1,\delta_2>0$. The event rule \cite[eq.~(18)]{Yang2025} is equivalently
\begin{equation}\label{eq:trigger}
 t_{i+1}=\min\bigl\{t\in\{n,\ldots,L\}:D_m(t)\ge\Theta(t)\bigr\}.
\end{equation}
Here the underscore in $\zm$ and $\zp$ denotes the compact output block, as in \cite{Yang2025}.

\begin{proposition}\label{prop:empty}
Suppose the initial state and prior input-output history are zero and the actual measurement noise vanishes. For the predictive optimization \cite[eq.~(17)]{Yang2025}, the zero candidate $(\bar u,\underline{\bar z},\alpha,\sigma)=0$ is feasible and optimal: its objective value is zero, and all cost terms are nonnegative. For any $\epsb>0$ and $\delta_1>0$, the event set in \eqref{eq:trigger} is empty along the resulting zero trajectory.
\end{proposition}
\begin{proof}
The predicted and measured output blocks both vanish at every candidate time. The event inequality consequently reduces to $0\ge\delta_1\epsb$, which is impossible.
\end{proof}

The finite predicted control sequence can thus be exhausted without another trigger being scheduled by \eqref{eq:trigger}. A complete event implementation may use the terminal convention $t_{i+1}=L$ if no earlier crossing occurs; such a convention also defines a terminal case that must be represented in the closed-loop analysis. The example concerns well-posed scheduling under the displayed rule and makes no assertion of physical instability.

We next consider the critical step in the published UUB proof. At a genuine triggering instant, \eqref{eq:trigger} gives a \emph{lower} bound for $D_m(t_{i+1})$. Nevertheless, \cite[eq.~(26)]{Yang2025} asserts an \emph{upper} bound
\begin{equation}\label{eq:claimed}
 \begin{split}
 D_n(t_{i+1})&\le\delta_1\epsb
 +\delta_2\|\zp_{t_{i+1}-1}\|^2
 +(r+1)t_{i+1}\epsb^2,\\
 D_n(t)&:=\|\zn_{[0,t-1]}\|^2-\|\zp_{[0,t-1]}\|^2.
 \end{split}
\end{equation}
The relation $\zm=\zn+\varepsilon$ and the noise bound do not reverse the inequality in \eqref{eq:trigger}.

With zero measurement noise, $D_n=D_m$. At the first eligible index $t=n$, the trigger gives only $D_m(n)\ge\Theta(n)$. Set the predicted block to zero and the measured prefix norm to $M$. For sufficiently large $M$, the trigger holds, whereas \eqref{eq:claimed} would require
\[
 M^2\le\delta_1\epsb+(r+1)n\epsb^2.
\]
Thus the trigger and noise bound alone do not imply \eqref{eq:claimed}; a separate threshold-crossing estimate is needed.

More generally, if $t_{i+1}>n$ is the first triggering index and
\[
 a(t)=\|\zm_{t-1}\|^2-\|\zp_{t-1}\|^2,
\]
then $D_m(t)=D_m(t-1)+a(t)$. The absence of a trigger at $t-1$ provides only
\begin{equation}\label{eq:overshoot}
 D_m(t)<\Theta(t-1)+a(t),
 \qquad t=t_{i+1}>n.
\end{equation}
The threshold at $t-1$ depends on the measured and predicted prefix only through that earlier time; its comparison with $D_m(t)$ involves the newly included block at $t-1$. A usable upper estimate at the event instant therefore requires control of the one-step crossing increment $a(t)$, together with the noise terms. Neither this overshoot estimate nor an analysis of an explicitly forced terminal event is supplied in the passage from (18) to (26) in \cite{Yang2025}. The subsequent recursion (27)--(32) consequently does not establish the asserted UUB bound from the displayed hypotheses.

\section{Conclusions}
The realizable geometry of multirate output blocks yields a universal rank obstruction to square inverse-based compact lifting and an exact reconstruction through rank-complete data. The full family of predictor extensions admits a canonical Moore--Penrose member with minimum Frobenius norm and minimum worst-case Euclidean measurement-noise amplification. The joint data predictor is invariant across rank-complete exact experiments and can be identified from mixed input--output records with arbitrary initial states. Quantitative bounds track identification errors through informative singular-value margins, while the reduced predictor preserves the nominal quadratic cost. Complementary results locate the relevant observability test at the fast sampling rate and the threshold-crossing estimate required in discrete event-triggered stability analysis. Together these statements provide a geometric, constructive, and perturbation-aware framework for compact multirate prediction.

\end{document}